# A Geometric Approach to Solve Fuzzy Linear Systems of Differential Equations

N. Gasilov [a], Sh. G. Amrahov (Şahin Emrah) [b], A. Golayoğlu Fatullayev [c]

[a] Baskent University, Ankara, 06810 Turkey (phone: (+90)3122341010/1220; fax: (+90)3122341051; e-mail: gasilov@baskent.edu.tr). (corresponding author)

[b] Department of Computer Engineering, Engineering Faculty of Ankara University, Aziz Kansu Building, Tandogan Kampus, Ankara, 06100 Turkey (e-mail: emrah@eng.ankara.edu.tr).

[c] Baskent University, Ankara, 06810 Turkey (e-mail: afet@baskent.edu.tr).

In this paper, systems of linear differential equations with crisp real coefficients and with initial condition described by a vector of fuzzy numbers are studied. A new method based on the geometric representations of linear transformations is proposed to find a solution. The most important difference between this method and methods offered in previous papers is that the solution is considered to be a fuzzy set of real vector-functions rather than a fuzzy vector-function. Each member of the set satisfies the given system with a certain possibility. It is shown that at any time the solution constitutes a fuzzy region in the coordinate space, $\alpha$-cuts of which are nested parallelepipeds. Proposed method is illustrated on examples.

## 1. Introduction

It is natural to model dynamical systems with uncertainty by fuzzy systems of differential equations (FSDE). Linear FSDE, in particular, appear in many applications.

The concept of a fuzzy derivative was defined by Chang and Zadeh in [11]. It was followed up by Dubois and Prade in [14], who used the extension principle. The term "fuzzy differential equation" was introduced in 1987 by Kandel and Byatt [22-23].

There have been many suggestions for the definition of fuzzy derivative to study fuzzy differential equations. One of the earliest suggestions was to generalize the Hukuhara derivative [17] of a set-valued function. This generalization was made by Puri and Ralescu [26] and studied by Kaleva [19, 20]. It soon appeared that the solution of a fuzzy differential equation defined by means of Hukuhara derivative has a deficiency: it becomes fuzzier as time goes, according to [13, 21]. Hence, behavior of the fuzzy solution is quite different from that of the crisp solution. Seikkala [28] introduced the notion of fuzzy derivative as an extension of the Hukuhara derivative and fuzzy integral, which was the same as what Dubois and Prade [14] proposed. To circumvent problems arising in connection with Hukuhara derivative, Hüllermeier [18] considered fuzzy differential equation as a family of differential inclusions. The main downside of using differential inclusions is that we do not have an adequate definition for derivative of a fuzzy-number-valued function. The concept of strongly generalized differentiability was introduced in [4] and studied in [5, 6, 9]. In [24] a generalized concept of higher-order differentiability for fuzzy functions is presented to solve *N*th-order fuzzy differential equations. Buckley and Feuring [7] and Buckley et al. [8] gave a very general formulation of fuzzy first-order initial value problem. They firstly find the crisp solution, fuzzify it and then check to see if it satisfies the FSDE. In [27] Rodriguez-Lopez considered several comparison results for the solutions of FSDE obtained through different methods using the Hukuhara derivative. Allahviranloo, Kiani and Motamedi [1] applied differential transformation method by using generalized H-differentiability. In [25] Mizukoshi et al. showed that the solutions of the Cauchy problem obtained by the Zadeh extension principle and by using a family of differential inclusions are same, but in [2] Allahviranloo,

Shafiee and Nejatbakhsh demonstrated with an example that the main result of [25] was incorrect. Xu, Liao and Hu in [29] used complex number representation of $\alpha$-level sets of a fuzzy system and proved theorems that provide the solutions in this representation. Chalco-Cano and Román-Flores [10] studied the class of fuzzy differential equations where the dynamics is given by a continuous fuzzy mapping which is obtained via Zadeh's extension principle.

In this paper, we apply a geometric approach to fuzzy linear system of differential equations (FLSDE) with crisp real coefficients and with initial condition described by a vector of fuzzy numbers. We interpret a vector of fuzzy numbers as a rectangular prism in $n$ dimensional space, and show that at any time the solution corresponds to an $n$ dimensional parallelepiped. Unlike earlier research, we are not looking for solutions of FLSDE in the form of fuzzy vector-function. Instead, our solutions constitute a fuzzy set of real vector-functions. Each member in the solution set satisfies the system with a certain possibility.

In articles [15, 16], by using the same geometric approach we proposed an algorithm to solve linear systems of algebraic equations with crisp coefficients and with fuzzy numbers on the right-hand side. Here we adopt the approach to the fuzzy linear system of differential equations.

This paper is comprised of 6 sections including the Introduction. Preliminaries are given in Section 2. In Section 3, we define FLSDE. In Section 4, we apply a geometric approach to find the solution of FLSDE and present the main results. In Section 5, we solve samples of FLSDE. In Section 6, we summarize the results.

## 2. Preliminaries

We define a fuzzy number $\tilde{u}$ in parametric form [12].

**Definition 1**. A fuzzy number $\tilde{u}$ in parametric form is a pair $(\underline{u}, \overline{u})$ of functions $\underline{u}(r), \overline{u}(r)$, $0 \leq r \leq 1$, which satisfies the following requirements:

1. $\underline{u}(r)$ is a bounded monotonic nondecreasing left continuous function over [0,1]
2. $\overline{u}(r)$ is bounded monotonic nonincreasing left continuous function over [0,1]
3. $\underline{u}(r) \leq \overline{u}(r), \ 0 \leq r \leq 1$

The set of all these fuzzy numbers is denoted by $E^1$.

A popular fuzzy number is the triangular number $\tilde{u} = (a, c, b)$ with the membership function

$$\mu(x) = \begin{cases} \dfrac{x-a}{c-a}, & a \leq x \leq c \\ \dfrac{x-b}{c-b}, & c \leq x \leq b \end{cases}$$

where $c \neq a$ and $c \neq b$. For triangular numbers we have $\underline{u}(r) = a + (c-a)r$ and $\overline{u}(r) = b + (c-b)r$.

We will denote $\underline{\underline{u}} = a$ and $\overline{\overline{u}} = b$ to indicate the left and the right limits of $\tilde{u}$, respectively.

We can represent a crisp number $a$ by taking $\underline{u}(r) = \overline{u}(r) = a$, $0 \leq r \leq 1$.

For two arbitrary fuzzy numbers $\tilde{u}$ and $\tilde{v}$ the equality $\tilde{u} = \tilde{v}$ means that $\underline{u}(r) = \underline{v}(r)$ and $\overline{u}(r) = \overline{v}(r)$.

For two arbitrary fuzzy numbers $\tilde{u}$ and $\tilde{v}$ the following arithmetic operations are defined:

*a*) Addition: $\tilde{u} + \tilde{v} = (\underline{u}(r) + \underline{v}(r), \overline{u}(r) + \overline{v}(r))$

*b*) Multiplication by a real number $k$:
$$k\tilde{u} = \begin{cases} (k\underline{u}(r), k\overline{u}(r)), & k \geq 0 \\ (k\overline{u}(r), k\underline{u}(r)), & k < 0 \end{cases}$$

*c*) Subtraction: $\tilde{u} - \tilde{v} = \tilde{u} + (-1)\tilde{v}$

## 3. Fuzzy linear systems of differential equations

**Definition 2.** Let $a_{ij}$, $(1 \leq i, j \leq n)$ be crisp numbers, $f_i(t), 1 \leq i \leq n$ be given crisp functions and $\tilde{u}_i = (\underline{u}_i(r), \overline{u}_i(r)), 0 \leq r \leq 1$, $1 \leq i \leq n$ be fuzzy numbers. The linear system of differential equations

$$\begin{cases} x_1'(t) = a_{11}x_1 + a_{12}x_2 + \ldots + a_{1n}x_n + f_1(t) \\ x_2'(t) = a_{21}x_1 + a_{22}x_2 + \ldots + a_{2n}x_n + f_2(t) \\ \vdots \\ x_n'(t) = a_{n1}x_1 + a_{n2}x_2 + \ldots + a_{nn}x_n + f_n(t) \end{cases} \quad (1)$$

with the fuzzy initial condition

$$\begin{cases} x_1(t_0) = \tilde{u}_1 \\ x_2(t_0) = \tilde{u}_2 \\ \vdots \\ x_n(t_0) = \tilde{u}_n \end{cases} \quad (2)$$

is called a fuzzy linear system of differential equations (FLSDE).

One can rewrite the problem (1)-(2) as follows using matrix notation.
$$\begin{cases} X' = AX + F(t) \\ X(t_0) = \tilde{B} \end{cases} \quad (3)$$

where $A = [a_{ij}]$ is an $n \times n$ crisp matrix, $F(t) = (f_1(t), f_2(t), \ldots, f_n(t))^T$ is a crisp vector-function and $\tilde{B} = (\tilde{u}_1, \tilde{u}_2, \ldots, \tilde{u}_n)^T$ is a vector of fuzzy numbers.

If differential equations are considered to describe motion of a body, then fuzzy initial conditions indicate some uncertainty regarding the location of the body at time $t_0$.

## 4. The method of solution

Now we show how to find a solution for the problem (3) as a fuzzy set.
Without loss of generality, we put $t_0 = 0$.
Let us write the initial value vector as $\tilde{B} = \mathbf{b}_{cr} + \tilde{\mathbf{b}}$, where $\mathbf{b}_{cr}$ is a vector with possibility 1 and denotes the vertex (crisp part) of the fuzzy region $\tilde{B}$, and $\tilde{\mathbf{b}}$ denotes the uncertain part with vertex at the origin. It is easy to see that, a solution of the given system is of the form

$\widetilde{X}(t) = \mathbf{x}_{cr}(t) + \widetilde{\mathbf{x}}(t)$ (crisp solution + uncertainty). Here $\mathbf{x}_{cr}(t)$ is a solution of the non-homogeneous crisp problem

$$\begin{cases} X' = AX + F(t) \\ X(0) = \mathbf{b}_{cr} \end{cases};$$

while $\widetilde{\mathbf{x}}(t)$ is a solution of the homogeneous system with fuzzy initial condition

$$\begin{cases} X' = AX \\ X(0) = \widetilde{\mathbf{b}} \end{cases}.$$

In regard to motion of a body, one could interpret $\mathbf{x}_{cr}(t)$ as the main trajectory. It is possible to compute $\mathbf{x}_{cr}(t)$ by means of analytical or numerical methods. Hence, (3) is reduced to solving a homogeneous system with fuzzy initial conditions.

The basis of our method is summarized below:

We shall make use of the following facts about linear transformations [3]:
1. A linear transformation maps the origin (zero vector) to the origin (zero vector).
2. A linear transformation maps a pair of parallel straight lines to a pair of parallel straight lines (thus a pair of parallel faces to a pair of parallel faces). Consequently, a linear transformation maps a parallelepiped to a parallelepiped.

In addition, we shall reference a property of fuzzy number vectors.
3. The fuzzy vector $\widetilde{\mathbf{b}}$ forms a fuzzy region in $R^n$, vertex of which is located at the origin and boundary of which is a rectangular prism. Furthermore, the $\alpha$-cuts of the region are rectangular prisms nested within one another.

The next three properties are in connection with the initial value problem

$$\begin{cases} \mathbf{x}' = A\mathbf{x} \\ \mathbf{x}(0) = \mathbf{v} \end{cases} \quad (4)$$

4. By the existence and uniqueness theorem, two solutions with two different initial points have distinct values for any $t$.
5. The solution of (4) is of the form $\mathbf{x}(t) = e^{At}\mathbf{v}$.
6. For a fixed $t = t_*$, $\mathbf{x}(t_*) = e^{At_*}\mathbf{x}(0) = M\mathbf{x}(0)$, where $M$ is a fixed invertible matrix. Hence, the value of the solution function at $t = t_*$ is determined by a linear transformation described by $M$.

The facts 1-6 bring us to the following conclusion: The set of initial points form a rectangular prism (or more generally speaking, a parallelepiped) and the points of the solution at any time form a parallelepiped.

In particular, for $n = 2$, rectangular prism and parallelepiped turn into rectangle and parallelogram, respectively. According to the discussion above, solution curves make up nested surfaces extending along $t$ direction (like a coaxial cable). Cross-sections of these surfaces at any $t = t_*$ are nested parallelograms.

We can liken the behavior of the solution in *xy* plane to that of a cloud of dust. In the cloud, there exists a point ("the center") where the density of dust is the highest. As moving away from the center, the density decreases along the parallelograms. In other words, the parallelograms correspond to the level curves of the density. The motion of the center is governed by the crisp solution. The cloud of dust moves along with the center; but the parameters of the parallelograms (orientations, ratio of sides) may change in time.

The foregoing discussion qualitatively indicates how the solution would behave in *n* dimensional coordinate space, in general. Now, we shall find a formula for the solution.

Firstly, we consider the case where the initial values in (2) are triangular fuzzy numbers $\tilde{u}_i = (l_i, m_i, r_i)$. We have $(b_{cr})_i = m_i$ and $\tilde{b}_i = (\underline{b}_i, 0, \overline{\overline{b}}_i) = (l_i - m_i, 0, r_i - m_i)$. Let us denote $\tilde{\mathbf{b}} = (\underline{\mathbf{b}}, \mathbf{0}, \overline{\overline{\mathbf{b}}})$, where $\underline{\mathbf{b}} = (\underline{b}_1, \underline{b}_2, \ldots, \underline{b}_n)$ and $\overline{\overline{\mathbf{b}}} = (\overline{\overline{b}}_1, \overline{\overline{b}}_2, \ldots, \overline{\overline{b}}_n)$.

One can express the vectors $\underline{\mathbf{b}}$ and $\overline{\overline{\mathbf{b}}}$ through standard basis vectors $\mathbf{e}_1, \mathbf{e}_2, \ldots, \mathbf{e}_n$:

$$\underline{\mathbf{b}} = \underline{b}_1 \mathbf{e}_1 + \underline{b}_2 \mathbf{e}_2 + \ldots + \underline{b}_n \mathbf{e}_n; \qquad \overline{\overline{\mathbf{b}}} = \overline{\overline{b}}_1 \mathbf{e}_1 + \overline{\overline{b}}_2 \mathbf{e}_2 + \ldots + \overline{\overline{b}}_n \mathbf{e}_n.$$

Let $\mathbf{v}_i = \underline{b}_i \mathbf{e}_i$ and $\mathbf{u}_i = \overline{\overline{b}}_i \mathbf{e}_i$. Note that $\mathbf{v}_i$ and $\mathbf{u}_i$ are vectors with all but *i*-th coordinates zero. The *i*-th coordinate of $\mathbf{v}_i$ is negative, while the *i*-th coordinate of $\mathbf{u}_i$ is positive. Any crisp vector in $R^n$ can be expressed uniquely as a linear combination, with positive coefficients, of the vectors $\mathbf{v}_i$ and $\mathbf{u}_i$.

The fuzzy initial vector $\tilde{\mathbf{b}}$ forms a rectangular prism in the coordinate space:

$$\tilde{\mathbf{b}} = \{\mathbf{b} = \alpha_1 \mathbf{w}_1 + \alpha_2 \mathbf{w}_2 + \ldots + \alpha_n \mathbf{w}_n \mid \alpha_i \in [0, 1]; \mathbf{w}_i = \mathbf{v}_i \text{ or } \mathbf{w}_i = \mathbf{u}_i\}$$

with membership function $\mu_{\tilde{\mathbf{b}}}(\mathbf{b}) = 1 - \max_{1 \leq i \leq n} \alpha_i$.

Let $\mathbf{q}_i(t) = e^{At} \mathbf{v}_i$ and $\mathbf{p}_i(t) = e^{At} \mathbf{u}_i$. Then the solution can be expressed as follows:

$$\tilde{X} = \{\mathbf{x}(t) = \mathbf{x}_{cr}(t) + \alpha_1 \mathbf{r}_1(t) + \alpha_2 \mathbf{r}_2(t) + \ldots + \alpha_n \mathbf{r}_n(t) \mid \alpha_i \in [0, 1]; \mathbf{r}_i = \mathbf{q}_i \text{ or } \mathbf{r}_i = \mathbf{p}_i\}; \qquad (5)$$

$$\mu_X(\mathbf{x}(t)) = 1 - \max_{1 \leq i \leq n} \alpha_i. \qquad (6)$$

If initial values are triangular fuzzy numbers, then one could determine an $\alpha$-cut of the solution, $X_\alpha$, through geometric similarity (with coefficient $1 - \alpha$) without additional computation.

In this paper, we define the possibility of a member of the solution set to be the possibility of the corresponding initial value.

Below we give another representation for the solution and then we generalize the results for the case where the initial values of the problem (3) are arbitrary fuzzy numbers.

The rectangular prism, corresponding to the initial values in the form of triangular fuzzy numbers, and its $\alpha$-cuts can also be represented as follows:

$$\tilde{\mathbf{b}} = \{c_1 \mathbf{e}_1 + c_2 \mathbf{e}_2 + \ldots + c_n \mathbf{e}_n \mid c_i \in [\underline{b}_i, \overline{\overline{b}}_i]\};$$

$$\mathbf{b}_\alpha = \{c_1 \mathbf{e}_1 + c_2 \mathbf{e}_2 + \ldots + c_n \mathbf{e}_n \mid c_i \in [(1-\alpha)\underline{b}_i, (1-\alpha)\overline{\overline{b}}_i]\}.$$

Let $\mathbf{g}_i(t) = e^{At}\mathbf{e}_i$. Then we can obtain the following formulas for $\alpha$-cuts of the solution and the solution itself:

$$X_\alpha = \{\mathbf{x}(t) = \mathbf{x}_{cr}(t) + c_1\mathbf{g}_1(t) + c_2\mathbf{g}_2(t) + \ldots + c_n\mathbf{g}_n(t) \mid c_i \in [(1-\alpha)\underline{b_i}, (1-\alpha)\overline{b_i}]\}; \quad (7)$$

$$\widetilde{X} = X_0 \text{ with } \mu_X(\mathbf{x}(t)) = 1 - \max_{1 \leq i \leq n} \gamma_i, \text{ where } \gamma_i = \begin{cases} c_i/\overline{b_i}, & c_i \geq 0 \\ c_i/\underline{b_i}, & c_i < 0 \end{cases}. \quad (8)$$

For the case when initial values consist of parametric fuzzy numbers $\tilde{u}_i = (\underline{u}_i(r), \overline{u}_i(r))$ the solution can be represented as follows:

$$X_\alpha = \{\mathbf{x}(t) = \mathbf{x}_{cr}(t) + c_1\mathbf{g}_1(t) + c_2\mathbf{g}_2(t) + \ldots + c_n\mathbf{g}_n(t)$$
$$\mid c_i \in [\underline{u}_i(\alpha) - (b_{cr})_i, \overline{u}_i(\alpha) - (b_{cr})_i]\} \quad (9)$$

$$\widetilde{X} = X_0 \quad (10)$$

with

$$\mu_X(\mathbf{x}(t)) = \min_{1 \leq i \leq n} \alpha_i, \text{ where } \alpha_i = \begin{cases} \overline{u}_i^{-1}(k_i), & k_i > \overline{u}_i(1)] \\ 1, & \underline{u}_i(1) \leq k_i \leq \overline{u}_i(1) \text{ and } k_i = (b_{cr})_i + c_i. \\ \underline{u}_i^{-1}(k_i), & k_i < \underline{u}_i(1) \end{cases} \quad (11)$$

Hence, to determine the solution set we need to calculate $e^{At}$. Note that $\mathbf{g}_i(t) = e^{At}\mathbf{e}_i$ is the solution of the crisp homogeneous system with initial value vector $\mathbf{e}_i$. In other words, we need to determine $n$ different solutions of the homogeneous system by taking the basis vectors $\mathbf{e}_i$, $1 \leq i \leq n$ as initial value vectors.

To summarize, the solution of the problem (1)-(2) (or, (3) in matrix form) is the fuzzy set of real vector-functions (or, fuzzy bunch of vector-functions), which can be represented by the formulas (9)-(11) in general case. If triangular fuzzy numbers describe initial conditions, the formulas (7)-(8) can be applied. To determine the fuzzy solution set we only need to work out the solution of crisp initial value problem for non-homogeneous system and to calculate $e^{At}$, or to find $n$ solution $\mathbf{g}_i(t)$, $(i = 1, \ldots, n)$ of the crisp homogeneous system.

We note that if the initial values are in parametric form, then $\mathbf{b}_{cr}$, in general, is not unique. In this case, we can choose the components of $\mathbf{b}_{cr}$ arbitrarily to the extent that $\underline{u}_i(1) \leq (b_{cr})_i \leq \overline{u}_i(1)$. For instance, we can put $(b_{cr})_i = [\underline{u}_i(1) + \overline{u}_i(1)]/2$.

## 5. Examples

**Example 1.**

Solve the system $\begin{bmatrix} x' \\ y' \end{bmatrix} = \begin{bmatrix} 3 & -1 \\ 4 & -2 \end{bmatrix} \begin{bmatrix} x \\ y \end{bmatrix} + \begin{bmatrix} 5t^2 - 15t - 25 \\ 10t^2 - 10t - 40 \end{bmatrix}$ (12)

with the initial values $\begin{bmatrix} x(0) \\ y(0) \end{bmatrix} = \begin{bmatrix} (14.5, 15, 16) \\ (4,\ 6,\ 9) \end{bmatrix}$.

Find the $\alpha$-cut of the solution for $\alpha = 0.5$.

**Solution.**
Initial values, given by triangular fuzzy numbers, form a fuzzy region in coordinate space. The boundary of this region (the rectangle *ABCD*) and $\alpha$-cut for $\alpha = 0.5$ is shown in Figure 1. The vertex of the region, $\mathbf{b}_{cr} = (15,\ 6)^T$, is marked with a dot.

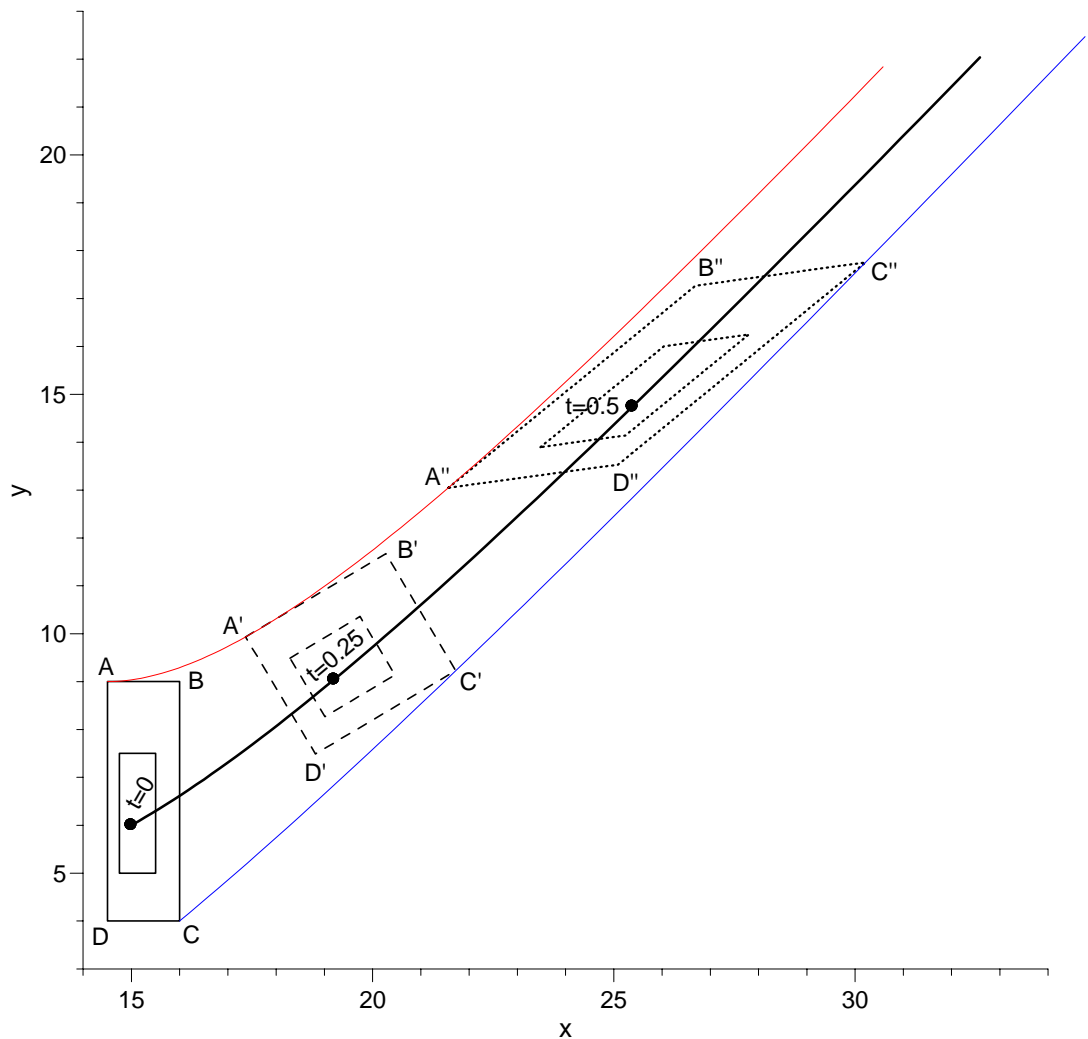

**Figure 1**. *ABCD* is the boundary of the fuzzy region corresponding to the initial values. Dashed parallelograms are the boundaries of the fuzzy region, corresponding to the solution, and its $\alpha = 0.5$-cut at time $t = 0.25$. Dotted parallelograms are the same, but for time $t = 0.5$. The thick line is the crisp solution.

The problem could be solved in two steps:

**1.** We determine the crisp solution corresponding to the non-homogenous system with the crisp initial values.

The solution of (12) with the initial value $\begin{bmatrix} x(0) \\ y(0) \end{bmatrix} = \begin{bmatrix} 15 \\ 6 \end{bmatrix}$ is $\begin{aligned} x_{cr}(t) &= 5(t+2) + \tfrac{1}{3}e^{-t} + \tfrac{14}{3}e^{2t} \\ y_{cr}(t) &= 5t^2 + \tfrac{4}{3}e^{-t} + \tfrac{14}{3}e^{2t} \end{aligned}$

The solution is graphed with the thick line in Figure 1.

**2.** We look for the solution corresponding to the homogeneous system with fuzzy initial value, the vertex of which is at the origin.
This means that we are looking for the fuzzy solution of

$$\begin{bmatrix} x' \\ y' \end{bmatrix} = \begin{bmatrix} 3 & -1 \\ 4 & -2 \end{bmatrix} \begin{bmatrix} x \\ y \end{bmatrix} \tag{13}$$

with initial value

$$\begin{bmatrix} x(0) \\ y(0) \end{bmatrix} = \begin{bmatrix} (-0.5, 0, 1) \\ (-2, 0, 3) \end{bmatrix} \tag{14}.$$

The rectangle $A_0 B_0 C_0 D_0$ determined by (14) can be obtained by translation of $ABCD$ by the vector $-\mathbf{b}_{cr}$. The vertices are:
$A_0(-0.5, 3), B_0(1, 3), C_0(1, -2), D_0(-0.5, -2)$
The general solution of (13) is:
$$\begin{cases} x(t) = c_1 e^{-t} + c_2 e^{2t} \\ y(t) = 4c_1 e^{-t} + c_2 e^{2t} \end{cases}.$$
For some initial point $P(a, b)$, the constants have the following values:
$c_1 = (-a+b)/3$ and $c_2 = (4a-b)/3$

Based on the solution of (13), one could work out $A_0', B_0', C_0'$, the locations of $A_0, B_0, C_0$ at time $t$, and hence obtain a parallelogram with three of its vertices at $A_0', B_0', C_0'$. If this parallelogram is translated by $\mathbf{x}_{cr}(t)$, we wind up with the boundary of region $(A'B'C'D')$ determined by the fuzzy solution at time $t$ (Figure 1).

Since the initial values are triangular numbers, one can find any $\alpha$-cut of the solution by geometric similarity without having to do additional computations (Figure 1). The formula for $\alpha$-cuts of the solution set $\tilde{X}$ can be obtained from (7):

$$X_\alpha = \left\{ \begin{bmatrix} x(t) \\ y(t) \end{bmatrix} = \begin{bmatrix} (15(t+2) + (1-c_1+c_2)e^{-t} + (14+4c_1-c_2)e^{2t})/3 \\ (15t^2 + 4(1-c_1+c_2)e^{-t} + (14+4c_1-c_2)e^{2t})/3 \end{bmatrix} \right.$$
$$\left. \bigg| -0.5(1-\alpha) \le c_1 \le (1-\alpha), \ -2(1-\alpha) \le c_2 \le 3(1-\alpha) \right\}$$

**Example 2.**
Solve the initial value problem and find the $\alpha$-cut of the solution for $\alpha = 0.5$.

$$\begin{bmatrix} x' \\ y' \end{bmatrix} = \begin{bmatrix} 3 & -1 \\ 4 & -2 \end{bmatrix} \begin{bmatrix} x \\ y \end{bmatrix} + \begin{bmatrix} 5t^2 - 15t - 25 \\ 10t^2 - 10t - 40 \end{bmatrix}$$

$$\begin{bmatrix} x(0) \\ y(0) \end{bmatrix} = \begin{bmatrix} (14.5 + 0.2r,\ 16 - 0.6r^2) \\ (4 + 1.75r^2,\ 9 - 2.5\sqrt{r}) \end{bmatrix} = \begin{bmatrix} (\underline{a}(r),\ \overline{a}(r)) \\ (\underline{b}(r),\ \overline{b}(r)) \end{bmatrix}$$

**Solution.**
We note that the given system is same as in Example 1, but with different initial values, which are represented in parametric form.

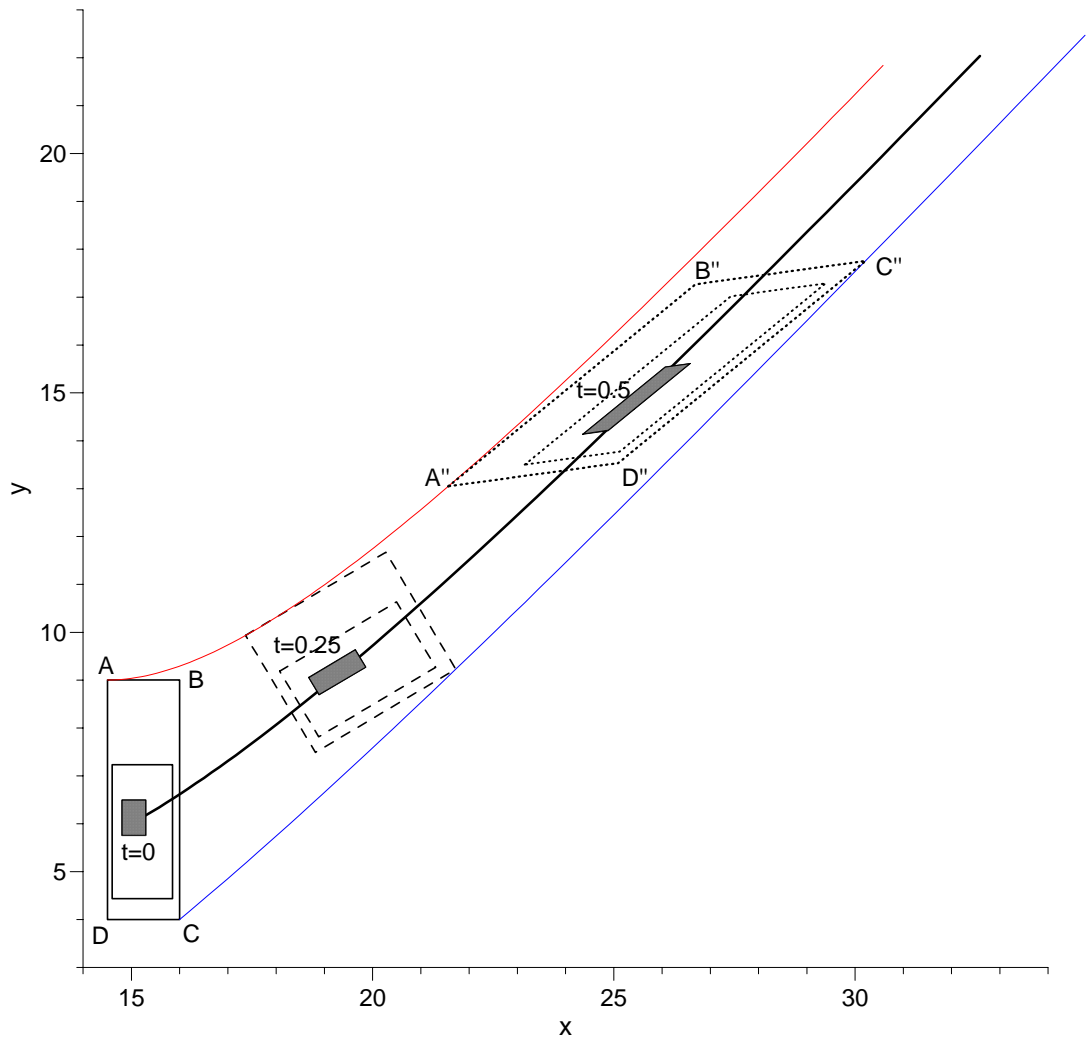

**Figure 2**. *ABCD* is the boundary of the fuzzy region corresponding to the initial values. Dashed parallelograms are the boundaries of the fuzzy region, corresponding to the solution, and its $\alpha = 0.5$-cut at time $t = 0.25$. Dotted parallelograms are the same, but for time $t = 0.5$. The shaded regions denote $\alpha = 1$ cuts. The thick line is the crisp solution.

For $r=1$ we have $\begin{bmatrix}(\underline{a}(1),\ \overline{a}(1))\\(\underline{b}(1),\ \overline{b}(1))\end{bmatrix}=\begin{bmatrix}(14.7,\ 15.4)\\(5.75,\ 6.50)\end{bmatrix}$. Therefore, we can arbitrarily choose the first and the second components of $\mathbf{b}_{cr}$ from the intervals $[14.7, 15.4]$ and $[5.75, 6.50]$, respectively. In the computations, we take $\mathbf{b}_{cr}=(15,\ 6)^T$, as in Example 1. Hence, formulas found in Example 1 are still valid.

To determine $\alpha$-cut of the solution, firstly we need to calculate the lower and upper limits of the corresponding $\alpha$-cut of the initial values. For $r=0.5$ we have
$\begin{bmatrix}(\underline{a}(0.5),\ \overline{a}(0.5))\\(\underline{b}(0.5),\ \overline{b}(0.5))\end{bmatrix}\approx\begin{bmatrix}(14.600,\ 15.8500)\\(4.4375,\ 7.23223)\end{bmatrix}$.

The results of the computation are shown in Figure 2. The values of the solution on the shaded regions have the possibility 1.

In Figure 2, the parallelograms, which are boundaries of $\alpha$-cuts, are not similar to the parallelogram corresponding to the solution. This is because the initial values are in parametric form rather than triangular fuzzy numbers.

The formula for $\alpha$-cuts of the solution set $\tilde{X}$ can be obtained from (9):

$$X_\alpha = \left\{\begin{bmatrix}x(t)\\y(t)\end{bmatrix}=\begin{bmatrix}(15(t+2)+(1-c_1+c_2)e^{-t}+(14+4c_1-c_2)e^{2t})/3\\(15t^2+4(1-c_1+c_2)e^{-t}+(14+4c_1-c_2)e^{2t})/3\end{bmatrix}\right.$$

$$\left.\vphantom{\begin{bmatrix}x\\y\end{bmatrix}}\ 0.2\alpha-0.5\leq c_1\leq 1-0.6\alpha^2,\ 1.75\alpha^2-2\leq c_2\leq 3-2.5\sqrt{\alpha}\right\}$$

The following example is taken from [8].

**Example 3. Arms race model**
Solve the initial value problem:
$$\begin{bmatrix}x'\\y'\end{bmatrix}=\begin{bmatrix}-3 & 2\\3 & -4\end{bmatrix}\begin{bmatrix}x\\y\end{bmatrix}+\begin{bmatrix}1\\2\end{bmatrix}$$
$$\begin{bmatrix}x(0)\\y(0)\end{bmatrix}=\begin{bmatrix}(70,100,130)\\(70,100,130)\end{bmatrix}$$

**Solution:**
We first determine the crisp solution. Consider the following:

$$\begin{bmatrix}x'\\y'\end{bmatrix}=\begin{bmatrix}-3 & 2\\3 & -4\end{bmatrix}\begin{bmatrix}x\\y\end{bmatrix}+\begin{bmatrix}1\\2\end{bmatrix}$$
$$\begin{bmatrix}x(0)\\y(0)\end{bmatrix}=\begin{bmatrix}100\\100\end{bmatrix}$$

The solution is:
$x_{cr}(t)=\frac{4}{3}+98\frac{3}{5}e^{-t}+\frac{1}{15}e^{-6t}$
$y_{cr}(t)=\frac{3}{2}+98\frac{3}{5}e^{-t}-\frac{1}{10}e^{-6t}$

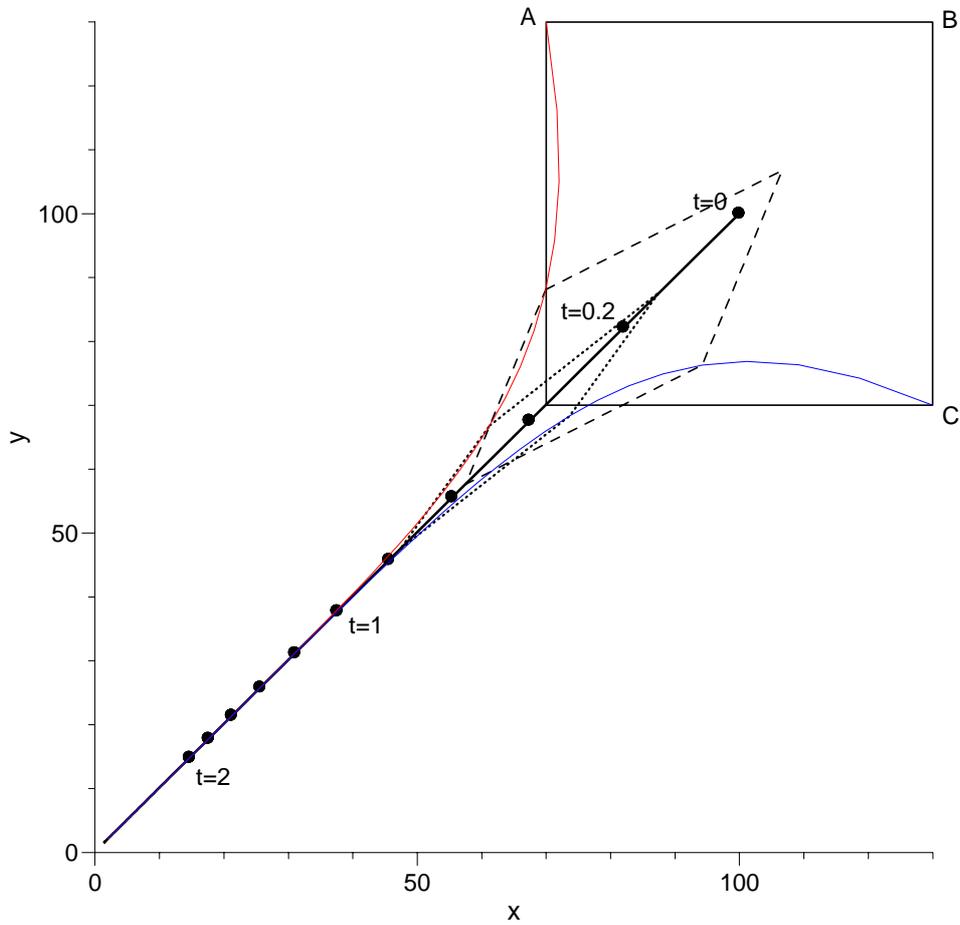

**Figure 3**. Rectangle *ABCD* is the boundary of the fuzzy region corresponding to the initial values. Dashed parallelogram is the boundary of the fuzzy region, corresponding to the solution at time $t = 0.2$. Dotted parallelogram is the same, but for time $t = 0.4$. The thick line is the crisp solution.

Secondly, we consider:

$$\begin{bmatrix} x' \\ y' \end{bmatrix} = \begin{bmatrix} -3 & 2 \\ 3 & -4 \end{bmatrix} \begin{bmatrix} x \\ y \end{bmatrix} + \begin{bmatrix} 1 \\ 2 \end{bmatrix}$$

$$\begin{bmatrix} x(0) \\ y(0) \end{bmatrix} = \begin{bmatrix} (-30, 0, 30) \\ (-30, 0, 30) \end{bmatrix}$$

We use the notation introduced in Example 1. Then,

$A_0(-30, 30)$, $B_0(30, 30)$, $C_0(30, -30)$, $D_0(-30, -30)$

The solution of the homogeneous system corresponding to initial point $P(a, b)$ is given by:

$$\begin{cases} x_P(t) = c_1 e^{-t} + 2c_2 e^{-6t} \\ y_P(t) = c_1 e^{-t} - 3c_2 e^{-6t} \end{cases}$$

where $c_1 = (3a+2b)/5$ and $c_2 = (a-b)/5$.

As shown in Figure 3, the fuzzy region corresponding to the solution at time $t$ gets smaller as $t$ increases, and shrinks down to a point as $t$ approaches infinity.

The $\alpha$-cuts of the solution set $\widetilde{X}$ can be expressed by the following formula:

$$X_\alpha = \left\{ \begin{bmatrix} x(t) \\ y(t) \end{bmatrix} = \begin{bmatrix} \frac{4}{3} + (98.6 + 0.6c_1 + 0.4c_2)e^{-t} + (\frac{1}{15} + 0.4c_1 - 0.4c_2)e^{-6t} \\ \frac{3}{2} + (98.6 + 0.6c_1 + 0.4c_2)e^{-t} + (-\frac{1}{10} - 0.6c_1 + 0.6c_2)e^{-6t} \end{bmatrix} \right.$$
$$\left. \Big| -30(1-\alpha) \le c_1 \le 30(1-\alpha),\ -30(1-\alpha) \le c_2 \le 30(1-\alpha) \right\}$$

## 6. Conclusion

In this paper, we dealt with systems of linear differential equations with crisp real coefficients and fuzzy initial condition. We proposed a geometric approach to solve the problem. Instead of looking for the solution as a fuzzy vector-function, we determined the solution as a fuzzy set of vector-functions, each of which satisfies FLSDE with some possibility. We showed that at a given time the values of the solutions form an $n$ dimensional parallelepiped. We suggested an efficient method to compute the fuzzy solution set. We illustrated the results with numerical examples.

**References**

[1] T. Allahviranloo, N.A. Kiani, N. Motamedi, Solving fuzzy differential equations by differential transformation method, *Information Sciences*, 179 (2009) 956–966.

[2] T. Allahviranloo, M. Shafiee, Y. Nejatbakhsh, A note on ''Fuzzy differential equations and the extension principle'', *Information Sciences*, 179 (2009) 2049–2051.

[3] Anton, H., Rorres, C., Elementary Linear Algebra, Applications Version: 9th Edition, John Wiley & Sons, 2005.

[4] B. Bede, S.G. Gal, Almost periodic fuzzy-number-valued functions, *Fuzzy Sets and Systems*, 147 (2004) 385-403.

[5] B. Bede, S.G. Gal, Generalizations of the differentiability of fuzzy number value functions with applications to fuzzy differential equations, *Fuzzy Sets and Systems*, 151 (2005) 581-99.

[6] B. Bede, I.J. Rudas, A.L. Bencsik, First order linear fuzzy differential equations under generalized differentiability, *Information Sciences*, 177 (2007) 1648-1662.

[7] J.J. Buckley, T. Feuring, Fuzzy differential equations, *Fuzzy Sets and Systems*, 110 (2000) 43-54.

[8] J.J. Buckley, T. Feuring and Y. Hayashi, Linear Systems of First Order Ordinary Differential Equations: Fuzzy Initial Conditions, *Soft Computing*, 6 (2002) 415-421.

[9] Y. Chalco-Cano, H. Román-Flores, On new solutions of fuzzy differential equations, Chaos, *Solitons & Fractals*, 38 (2008) 112-119.

[10] Y. Chalco-Cano, H. Román-Flores, Comparison between some approaches to solve fuzzy differential equations, *Fuzzy Sets and Systems*, 160 (2009) 1517–1527.

[11] S.L. Chang, L.A. Zadeh, On fuzzy mapping and control, *IEEE Transactions on Systems Man Cybernetics*, 2 (1972) 330-340.

[12] W. Cong-Xin and M. Ming. "Embedding problem on fuzzy number space: Part III", *Fuzzy Sets and Systems*, 46 (1992) 281-286.

[13] P. Diamond, P. Kloeden, Metric Spaces of Fuzzy Sets, World Scientific, Singapore 1994.

[14] D. Dubois and H. Prade, Towards Fuzzy Differential Calculus: Part 3, Differentiation, *Fuzzy Sets and Systems*, 8 (1982) 225-233.

[15] N.A. Gasilov, Sh.G. Amrahov, A.G. Fatullayev, H.I. Karakaş, Ö. Akın, Existence theorem for fuzzy number solutions of fuzzy linear systems. Proceedings. pp. 436-439. 1st International Fuzzy Systems Symposium (FUZZYSS'09). TOBB University of Economics and Technology. October 1-2, 2009. Ankara, Turkey. Edited by T. Dereli, A. Baykasoglu and I.B. Turksen.

[16] N. Gasilov, Sh.G. Amrahov, A.G. Fatullayev, H.I. Karakaş, Ö. Akın, A geometric approach to solve fuzzy linear systems. Submitted to *Information Sciences*.

[17] M. Hukuhara, Intégration des applications measurables dont la valeur est un compact convexe, *Funkcial. Ekvac.*, 10 (1967) 205–223.

[18] E. Hüllermeier, An approach to modelling and simulation of uncertain systems, *Int. J. Uncertain. Fuzz., Knowledge-Based System*, 5 (1997) 117-137.

[19] O. Kaleva, Fuzzy differential equations, *Fuzzy Sets and Systems*, 24 (1987) 301-317.

[20] O. Kaleva, The Cauchy Problem for Fuzzy Differential Equations, *Fuzzy Sets Systems*, 35 (1990) 389-396.

[21] O. Kaleva, A note on fuzzy differential equations, *Nonlinear Anal.*, 64 (2006) 895-900.

[22] A. Kandel, W.J. Byatt, Fuzzy differential equations, in: *Proceedings of International Conference Cybernetics and Society*, Tokyo; 1978, pp. 1213-1216.

[23] A. Kandel, W.J. Byatt, Fuzzy processes, *Fuzzy Sets and Systems*, 4 (1980) 117-152.

[24] A. Khastan, F. Bahrami, K. Ivaz, New Results on Multiple Solutions for Nth-Order Fuzzy Differential Equations under Generalized Differentiability, *Boundary Value Problems*, 2009 (2009) 13 p. Article ID 395714, doi:10.1155/2009/395714.

[25] M.T. Mizukoshi, L.C. Barros, Y. Chalco-Cano, H. Román-Flores, R.C. Bassanezi, Fuzzy differential equations and the extension principle, *Information Sciences*, 177 (2007) 3627–3635.

[26] M.L. Puri, D.A. Ralescu, Differentials of fuzzy functions, *J. Math. Anal. Appl.*, 91 (1983) 552-558.

[27] R. Rodriguez-Lopez, Comparison results for fuzzy differential equations, *Information Sciences*, 178 (2008) 1756–1779.

[28] S. Seikkala, On the fuzzy initial value problem, *Fuzzy Sets and Systems*, 24(3) (1987) 319–330.

[29] J. Xu, Z. Liao, Z. Hu, A class of linear differential dynamical systems with fuzzy initial condition, *Fuzzy Sets and Systems*, 158(21) (2007) 2339-2358.